\documentclass[a4paper,final]{amsart} 
\usepackage{amsmath}
\usepackage{epsfig}
\usepackage{amssymb,showkeys} 
\usepackage{xspace}
\usepackage{cite}

\newcommand*{\Matlab}{\texttt{Matlab}\xspace}

\newcounter{mnotecount} %
\usepackage{xcolor}

\begin{document}

\title{Computational approach to hyperelliptic Riemann surfaces}

\author{J.~Frauendiener }
\email{joergf@maths.otago.ac.nz}
\address{Department of Mathematics and Statistics, 
University of Otago,      
P.O. Box 56, Dunedin 9010, New Zealand}
\author{C.~Klein}
    \email{Christian.Klein@u-bourgogne.fr}
\address{Institut de Math\'ematiques de Bourgogne,
		Universit\'e de Bourgogne, 9 avenue Alain Savary, 21078 Dijon
		Cedex, France}

\date{\today}    
\thanks{This work has been supported in part by the Marie-Curie IRSES 
project RIMMP and the Royal Society of New Zealand}

\begin{abstract}
  We present a computational approach to general hyperelliptic Riemann surfaces in
  Weierstrass normal form. The surface is either given by a list of the branch
  points, the coefficients of the defining polynomial or a system of 
  cuts for the curve. A
  canonical basis of the homology is introduced algorithmically for this curve. The
  periods of the holomorphic differentials and the Abel map are computed with the
  Clenshaw-Curtis method in order to achieve spectral accuracy. The code can handle
  almost degenerate Riemann surfaces. This work generalizes previous work on real
  hyperelliptic surfaces with prescribed cuts to arbitrary hyperelliptic surfaces.
  As an example, solutions to the sine-Gordon equation in terms of multi-dimensional
  theta functions are studied, also in the solitonic limit of these solutions.
\end{abstract}


\maketitle

\section{Introduction}
Among compact Riemann surfaces of higher genus, hyperelliptic surfaces have clearly
the most applications and can be seen as straight forward generalizations of
elliptic surfaces.  It can be shown for instance that all Riemann surfaces of genus
2 are hyperelliptic. For this reason alone and for their relative simplicity,
hyperelliptic Riemann surfaces are important examples in algebraic geometry. Their
first applications in mathematical physics were Neumann's \cite{neu} integration of
the motion of a rigid body on a surface and the Kowalewskaja top \cite{kowa}. In the
1970s, quasiperiodic solutions of the celebrated Korteweg-de Vries (KdV) equation
were given by Its and Matveev in terms of multidimensional theta functions associated to hyperelliptic
Riemann surfaces, see \cite{algebro} and \cite{Dintro} for an account of the
history. It could be shown that many $1+1$ dimensional completely integrable
equations such as KdV have \emph{finite gap} solutions defined on hyperelliptic
surfaces, see for instance \cite{algebro} and references therein for the 
nonlinear Schr\"odinger (NLS) and \cite{kotlyarov} for the 
sine-Gordon (SG) equation, \cite{ernstbook} and \cite{korot1} for the 
Ernst equation and \cite{kalla} and references therein for the 
Camassa-Holm equation. Hyperelliptic Riemann surfaces also appear in 
the asymptotic description of \emph{dispersive shocks}, i.e., highly 
oscillatory regions in the solution  for these 
$1+1$-dimensional equations, see e.g.\ \cite{LL,GT,kam2,GK} for KdV 
and NLS. In this paper, we present an efficient numerical approach to
general hyperelliptic Riemann surfaces even in almost degenerate situations.

As all compact Riemann surfaces, hyperelliptic surfaces can be defined via algebraic
curves which, in the present case, can be written in Weierstrass normal form
\begin{equation}
    \mu^{2}=(z-z_{1})(z-z_{2})\ldots (z-z_{N}),
    \label{mu}
\end{equation}
where $(\mu,z)\in \mathbb{C}^{2}$ and where $z_{i}\in \mathbb{C}$ (all $z_{i}$
distinct) for $i=1,\ldots,N$. For $N=1,2$, the surface has genus $g=0$, for $N=3,4$,
it is elliptic ($g=1$), for $N>4$ the surface is hyperelliptic and has the genus
$g=N/2-1$ for $N$ even and $g=(N-1)/2$ for $N$ odd.  Hyperelliptic surfaces in
Weierstrass form can thus be represented in a standard way as two-sheeted coverings
of the complex plane which are branched at the zeros $z_{i}$, $i=1,\ldots,N$ of the
polynomial on the right hand side of (\ref{mu}), the \emph{branch points} of the
surface. A general point on the two-sheeted covering is denoted by $(z,\mu)$. The
hyperelliptic involution $\sigma$ interchanges the sheets, $\sigma
(z,\mu)=(z,-\mu)$. The sheets of the covering will be denoted in the following as
$+$ and $-$ sheet where the sign of $\mu$ is defined at some base 
point $z_{B}$ with $\mu(z_{B})\neq0$.  In
the case of odd $N$, the surface is branched at infinity, otherwise all branch
points are finite.  Note that surfaces obtained from (\ref{mu}) via 
birational transformations of the form $(\mu,z)\to
(R_{1}(\mu,z),R_{2}(\mu,z))$, where $R_{1}(\mu,z)$ and $R_{2}(\mu,z)$ are rational
in both arguments, i.e., by replacing $\mu$ and $z$ in (\ref{mu}) by 
$R_{1}(\mu,z)$ respectively $R_{2}(\mu,z)$, obviously also represent a hyperelliptic algebraic curve. We will
here only discuss curves in the Weierstrass normal form (\ref{mu}) (for a discussion
of how to check whether a given algebraic curve is hyperelliptic and how to find the
transformation to Weierstrass form, see for instance the appendix of
\cite{RSbookdp}).

The many applications of hyperelliptic curves make efficient numerical procedures
necessary. Different approaches are known: treating the hyperelliptic curves as
special cases of general algebraic curves, see \cite{deco01,deconinck03} and
\cite{RSbookfk}; uniformization techniques as described 
in~\cite{BB,RSbook}; and in terms
of theta functions \cite{ER}. Hyperelliptic curves can be treated simply as special
cases of general algebraic curves, for instance with the approach in
\cite{RSbookfk}. However, the hyperellipticity implies two important simplifications
which make special approaches for hyperelliptic curves much more efficient: firstly, it
is not necessary to solve an algebraic equation in order to analytically continue the roots
$\mu$ in (\ref{mu}) which, as was shown in \cite{RSbookfk}, takes almost 80\% of
the computing time --- it is just a square root here. This can be efficiently computed
for several values at the same time, and the analytic continuation is simply the
removal of inappropriate sign changes on a system of algorithmically identified cuts on
the surface. Also, the homology basis can be fixed \emph{a priori}. Secondly, the holomorphic
differentials are known \emph{a priori} for hyperelliptic surfaces (see
for instance \cite{algebro}):  the one-forms
\begin{equation}
    \nu_{k}=\frac{z^{k-1}dz}{\mu(z)},\quad k=1,\ldots,g
    \label{nu}
\end{equation}
form a basis of the space of holomorphic one-forms on the surface.  Thus, the
differentials do not need to be determined by the algorithm.  This means that
the operations, which are most time consuming in the case of general algebraic curves, become
either almost trivial or unnecessary in the hyperelliptic case. Therefore, efficient numerical approaches for
hyperelliptic curves are almost 100 times faster than more general codes. Hence, given the
importance of the hyperelliptic case, it is well worth to consider optimized adapted
approaches.

In \cite{prd,cam,lmp}, we have presented a computational approach to real Riemann
surfaces with \emph{spectral convergence}, i.e., an exponential decrease of the
numerical error with the numerical resolution. This approach allowed us to study
modular properties of hyperelliptic functions. Loosely speaking, this is the
dependence of functions defined on a family of hyperelliptic surfaces on the branch
points $z_{i}$ of the defining equation (\ref{mu}), see e.g.~\cite{klkoko} and
references therein.  This is also important in the context of solutions to the Ernst
equation found by Korotkin \cite{korot1}, see \cite{ernstbook} for a review, where
branch points of the underlying surface are parametrized by the physical
coordinates. The codes in \cite{cam,lmp} are also able to address the case where
branch points almost coincide, a case known as the \emph{solitonic limit} in the
context of quasiperiodic solutions in terms of theta functions to integrable
equations (in this limit, the solutions are no longer periodic, but localized
multi-solitons).  Since the signs of the square root in these codes were only fixed
locally, they were limited to particular real Riemann surfaces, i.e., branch points that were
either all real or all pairwise complex conjugate. In addition, the branch points
had to be known \emph{a priori} and be arranged in pairs corresponding to the cuts of the
surface.

In the present paper, we present a robust implementation of the square root which
allows the treatment of general hyperelliptic curves. In addition, an algorithm is
given for selecting a basis of the homology for a given (unordered) list of branch
points.  The curve (\ref{mu}) can be given in one of 3 forms: i) as an algebraic
equation in the form
\begin{equation}
    \mu^{2}=a_{N}z^{N}+a_{N-1}z^{N-1}+\ldots+a_{0}
    \label{mu2}
\end{equation}
with $a_{i}$, $i=0,\ldots,N$ some complex constants, and the branch points are
determined numerically, ii) in the  form~(\ref{mu}), and iii) in
terms of the branch points of the curve (\ref{mu}) given as in \cite{cam,lmp} in
pairs. In the first two cases, we present an algorithm to obtain a system of cuts
for the surface. The only condition is that the branch points are sufficiently
separated (as detailed in the following sections). In the third case, they can
almost coincide which allows for the study of almost solitonic situations. 

As an
example we will study quasiperiodic solutions to the SG equation
\begin{equation}
    u_{\xi\eta} = 4\sin(u)
    \label{SG},
\end{equation}
where $\xi$, $\eta$, $u$ are real, and where the index denotes a partial 
derivative. The SG equation has been first discovered in the theory of 
surfaces  as the Gauss-Codazzi equation 
for surfaces of constant negative curvature, see 
\cite{bour,bianchi} or \cite{korobianchi,ernstbook} for a more modern 
presentation. Note that the complete integrability of the SG equation 
was already shown by Bianchi \cite{bianchi} who managed to introduce 
a spectral parameter into the Gauss-Weingarten system for the 
surface, which plays the role of a Lax pair. As illustrated for 
instance in \cite{ernstbook}, the surfaces with constant negative 
curvature corresponding to SG solitons include the pseudosphere and 
the helicoid. The SG equation was also studied as a quantum field theoretical model in~\cite{faddeev}.

The paper is organized as follows: In section~2, an algorithm is presented to
establish a system of cuts for a given list of branch points $z_{i}$,
$i=1,\ldots,N$. In section~3, the square root $\mu$ in (\ref{mu}) is numerically
computed and analytically continued along these cuts. Also, the periods of the holomorphic
differentials and the Abel map are computed with the Clenshaw-Curtis method. In
section~4, we study solutions to the SG equation in terms of multidimensional theta
functions on hyperelliptic curves. We add some concluding remarks in section~5.

\section{Construction of a basis of the homology}
\label{sec:align-branch-points}

In this section we explain how to construct a basis of the homology of the surface
if the hyperelliptic curve is given in the form (\ref{mu}) or (\ref{mu2}), i.e.,
essentially for a list of the branch points $z_{i}$, $i=1,\ldots,N$ in (\ref{mu}).

\subsection{Computation of the branch points}

If the hyperelliptic curve is given in the form (\ref{mu2}), i.e., via the
coefficients $a_{i}$, $i=0,\ldots,N$ of the polynomial in $z$, then the branch points are
exactly the points where $\mu$ vanishes, which means the zeros of the polynomial. A
standard way to compute the zeros of a polynomial numerically is to construct the
\emph{companion matrix}, a matrix which has this polynomial as its characteristic
polynomial. The eigenvalues of this matrix are by definition the zeros of the
original polynomial. In \Matlab, this approach is used via the command
\emph{roots}. Note that the eigenvalues are computed with machine precision
($10^{-16}$ in our case, but in practice limited to $\approx 10^{-14}$ due to
rounding errors), but this does not imply that the zeros of the polynomial are
obtained with machine precision. Problems arise if the zeros are too close to each
other or of higher multiplicity because of the finite precision used in the
computations, see for instance \cite{zeng} and references given therein for more
details.

To avoid such problems, we always assume that the problem is well conditioned, i.e.,
that all zeros have multiplicity 1 and are well enough separated. To characterize
the latter condition we consider the quotient $\Delta$ of the largest distance between any
two finite branch points divided by the smallest distance between any two branch
points,
\begin{equation}
    \Delta:= 
    \frac{\max_{i,j=1,\ldots,N}|z_{i}-z_{j}|}{\min_{i,j=1,\ldots,N; 
    i\neq j}|z_{i}-z_{j}|},
    \label{Delta}
\end{equation}
and require that $\Delta$ is always smaller than $10^{6}$. In this case the \Matlab
command \emph{roots} produces a list of the zeros of the polynomial and thus the
branch points. If considerably higher values of $\Delta$ are to 
be treated, the branch points will be determined with less and less 
accuracy. To obtain reliable results in this case, the
branch points have to be given in pairs as discussed below.

\subsection{Identifications of cuts and basis of the homology}

\begin{figure}[htb]
  \centering
  \includegraphics[width = 0.7\linewidth]{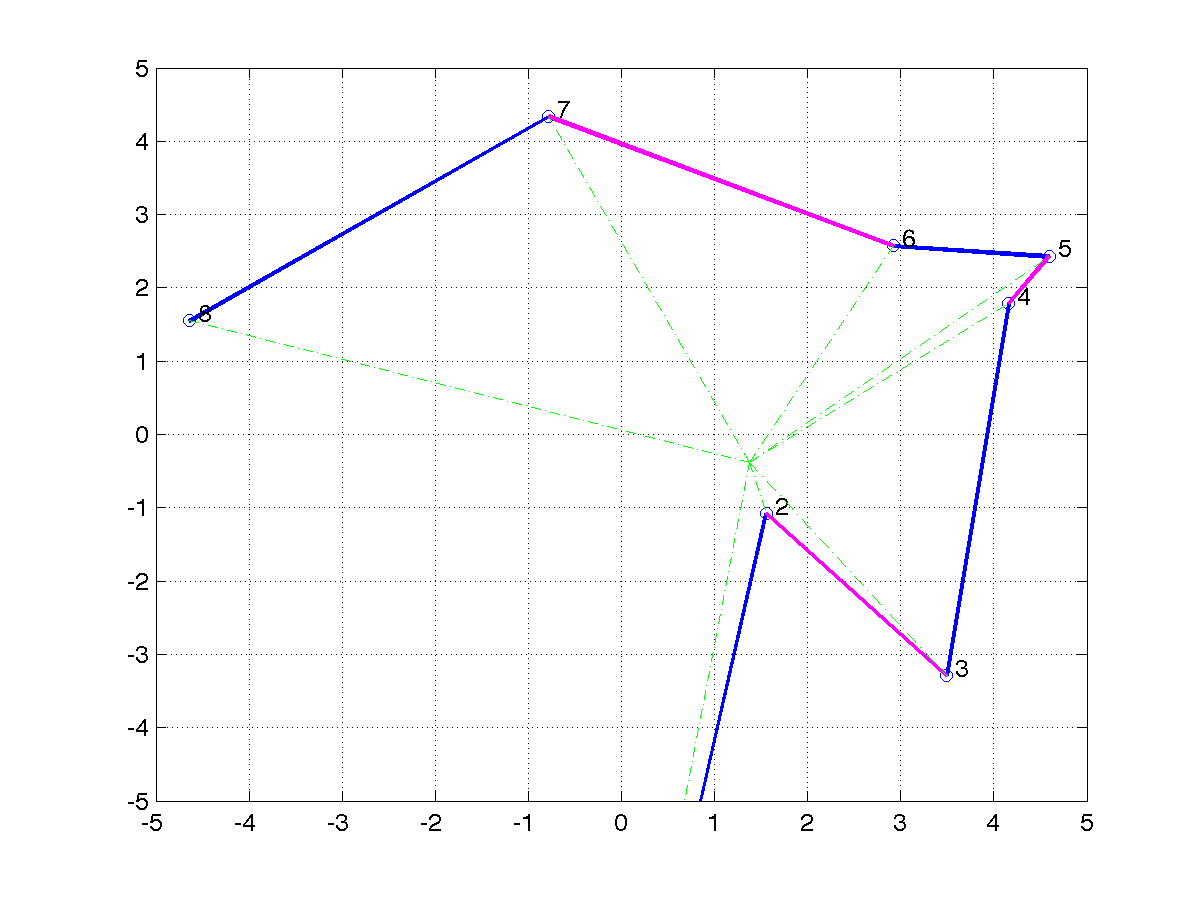}
  \caption{A generic situation. Six branch points chosen randomly are aligned using the argument with respect to their mean value. The blue lines indicate the cuts chosen on the basis of the sequence. The first branch point is at infinity.}
  \label{fig:genericcuts}
\end{figure}

Thus, whether the curve is given in the form (\ref{mu}) or (\ref{mu2}), we
essentially start in both cases with the set of branch points $z_i$. To define a
system of cuts between them, the branch points need to be sorted.  The general idea
is as follows: We define a function $f: \mathbb{C} \to \mathbb{R}$ which separates
the branch points, and then we sort them according to the values $f(z_i)$. This can
be done in several ways. Here, we restrict ourselves to the following 
method: We
determine a \emph{central point} $c$ taken as the arithmetic mean of the finite
branch points $z_{i}$, $i=1,\ldots,N$.  To reduce the possibility of degeneracies,
we add a random complex number of the order of machine precision to $c$.  Then we
compute the argument of the lines $z_{i}-c$, i.e., we put $f(z) = \arg(z-c)$.  We order the
branch points in increasing order of $f(z_i)$ into a sequence $(z_i)$ (in case of
degeneracies, the corresponding points are ordered according to their modulus, but
highly symmetric versions can be conveniently treated with the third approach
detailed below). Then we define the cuts $C_j$ to be between the branch points
$z_{2j+1}$ and $z_{2j+2}$ for $j=0,\ldots,g$, i.e., $C_0 = \overline{z_1z_2}$ etc. A
general picture is given in Fig.~\ref{fig:genericcuts}. Note, that the cuts tend to
form a line curving towards the left. If $N$ is odd, the surface is branched
at infinity, and we add $\infty$ as the first branch point by following the line through $c$ with
the argument of $z_{1}-c$ to infinity. This defines the first cut in this case.

\begin{figure}[htb]
  \centering
  \includegraphics[width=0.7\linewidth]{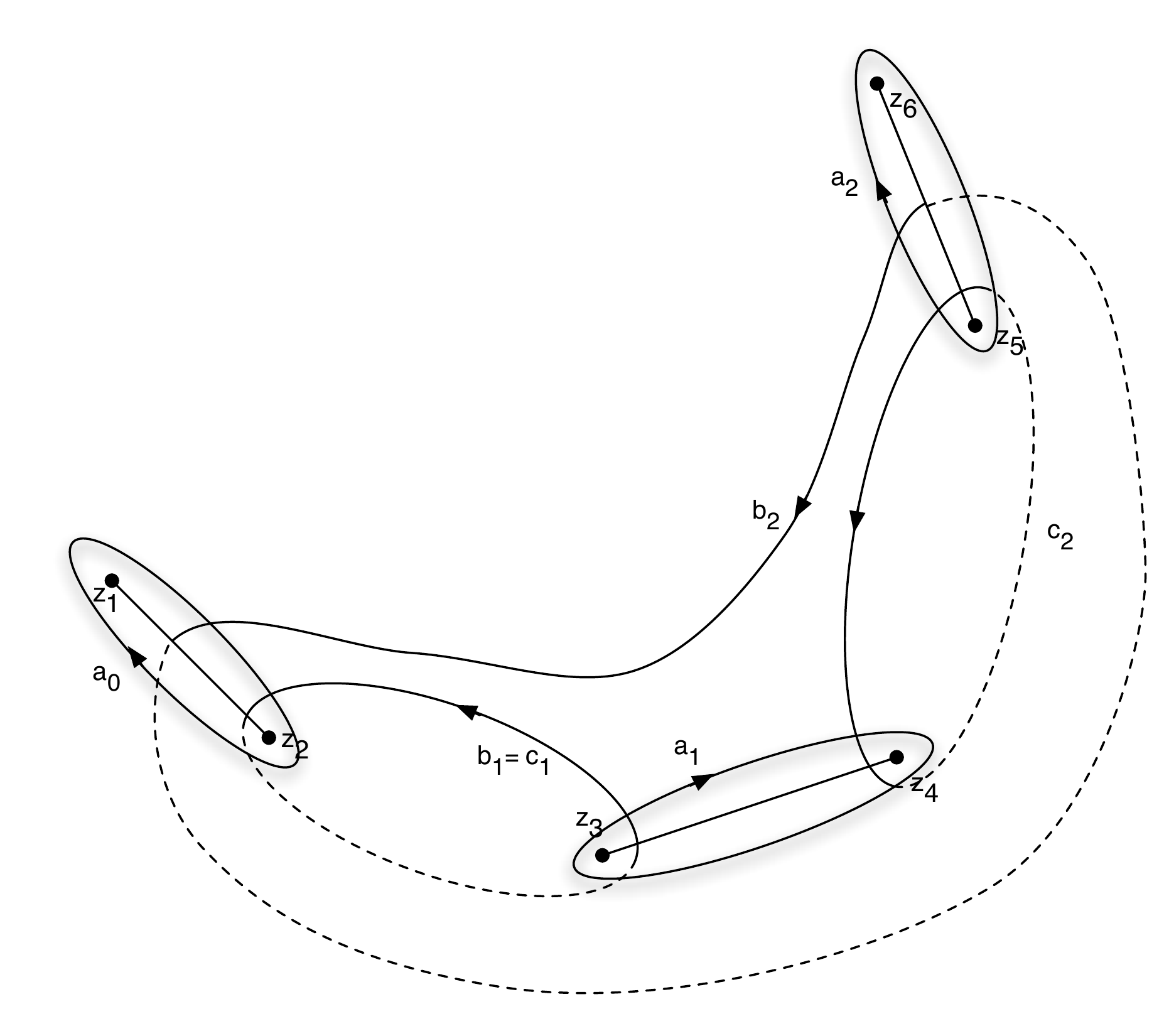}
  \caption{Basis of the homology used in the code  \label{fig:cutsystem}}
\end{figure}

To obtain a basis of the homology, the closed cycles are chosen according to the
cuts in Fig.~\ref{fig:genericcuts} as can be seen in Fig.~\ref{fig:cutsystem}. The
cycles denoted by $a_0$, $a_1$, \ldots $a_g$, enclose the cuts $C_j$, circling
around them exactly once in clock-wise orientation. We use $a_{1}, \ldots,a_{g}$ as
the so-called $a$-cycles.  Note that $a_{0}$ is homologous to the negative sum of the
$a_{i}$, $i=1,\ldots,g$ and thus not part of the canonical basis of the homology
used here. But it has a role in the ensuing computations and is, therefore, already
mentioned at this point.  The $b$-cycles $b_j$, $j=1,\ldots,g$, all start at the cut
$C_0$, extend on one sheet to the cut $C_j$ and then return on the other sheet to
$C_0$. Thus, $b_{j}$ intersects only the cycle $a_j$ and the intersection number
$a_{j}\circ b_{j}$ is $1$, when oriented in the indicated way. However, in the
computation of the periods in the following section, instead of dealing with these
cycles, we make use of the contours $c_j$, $j=1,\ldots,g$ which link successive cuts
as indicated in Fig.~\ref{fig:cutsystem}. Then one has the relation
\begin{equation}
b_j = \sum_{k=1}^{j} c_{k} .\label{eq:8}
\end{equation}

As mentioned in the introduction, the third option to use the code is to give the
branch points in the form of a list of pairs. These pairs will be interpreted by the
code as corresponding to the cuts $C_{j}$, $j=0,\ldots,g$ in
Fig.~\ref{fig:cutsystem} in ascending order. Thus with the choice of these pairs of
branch points, the system of cuts and thus the basis of the homology as in
Fig.~\ref{fig:cutsystem} are fixed. The cuts must be chosen in a way not to
intersect other cuts, otherwise the code will produce an error.

\section{Computation of Abel map and periods of the Riemann surface}

In this section we compute the periods of the hyperelliptic surface, i.e., the
integrals of the holomorphic one-forms along the cycles of the homology basis.  A
canonical basis $\omega_{k}$, $k=1,\ldots,g$ of holomorphic one-forms dual to a
canonical basis of the homology can be obtained from the basis (\ref{nu}) by the
normalization condition
\begin{equation}
    \int_{a_{j}}^{}\omega_{k}=2\pi \mathrm{i}\,\delta_{jk},\quad j,k=1,\ldots,g
    \label{normalization}.
\end{equation}

\subsection{Analytic continuation of the root}
\label{sec:cont-square-root}

To determine the periods of the surface, integrals of the form
\begin{equation}
  \label{eq:2}
  \int_\Gamma \nu_k = \int_\Gamma \frac{z^{k-1}}{\mu(z)}\,\mathrm{d} z
\end{equation}
for some contour $\Gamma$ on the hyperelliptic surface have to be computed. This
requires the analytic continuation of the root $\mu(z)$ defined in eq.~(\ref{mu})
along the contour. The general procedure for computing this function starts by
parametrizing the contour $\Gamma$ by a real parameter $s\in[-1,1]$. For
computational reasons, we choose a number $N_{c}+1$ of collocation points $s_{j}$,
$j=0,\ldots,N_{c}$ on the contour. We sample the function
\begin{equation}
    \tilde{\mu}(z):=\sqrt{\prod_{i=1}^{N}(z-z_{i})},\quad 
    \tilde{\mu}_{j}:= \tilde{\mu}(z(s_{j})),
    \label{mut}
\end{equation}
on the collocation points obtaining a vector $\tilde\mu_j$, where we use the \Matlab function \emph{sqrt} to determine the square root.  Note
that this square root is branched along the negative real axis, whereas $\mu(z)$ by
definition is an analytic function along the cuts defined in
Fig.~\ref{fig:genericcuts}. Thus, in general, the function $\tilde{\mu}(z)$ as
defined above will not be a smooth function on the contour $\Gamma$. 
To eliminate unwanted sign changes of $\tilde{\mu}(z)$, we determine 
the indices $j_{k}$ for which
\begin{equation}
  \label{eq:3}
  |\tilde{\mu}_{j+1}-\tilde{\mu}_{j}| > 
  |\tilde{\mu}_{j+1}+\tilde{\mu}_{j}|,
\end{equation}
i.e., the indices at which the \Matlab root is branched along the contour
$\Gamma$. This can be done in \Matlab efficiently in a vectorized way, which means
for all $j=1,\ldots,N_{c}$ at the same time. Then, for the $\tilde{\mu}_{j}$ with
$j_{k}\leq j\leq j_{k+1}$ and $k$ odd, the sign of the $\tilde{\mu}_{j}$ is
changed. The resulting function $\mu(z(s_{j}))$ gives the analytic continuation of
the root $\mu(z)$ along the contour $\Gamma$.  Since a sign change and thus a
quantity of order $2|\tilde{\mu}_{j}|$ has to be identified via (\ref{eq:3}), the
numerical result will be unique if the number $N_{c}$ of collocation points is
sufficiently large to reliably distinguish $\tilde{\mu}_{j}$ and $\pm
\tilde{\mu}_{j+1}$. This can be always achieved in practice as long as the branch
points are sufficiently separated, i.e., if $\Delta$ in (\ref{Delta}) 
is smaller than $10^{14}$.

\subsection{Computing the line integrals}
\label{sec:comp-line-integr}

To compute the periods of the holomorphic differentials $\nu_{k}$ in (\ref{nu}), the
latter have to be integrated over the cycles $a_{j}$ and $b_{j}$ in
Fig.~\ref{fig:cutsystem}.  However, it is preferable to integrate directly on the
cuts instead of closed contours around them in order to gain flexibility. This way,
the relative position of the branch points can be essentially arbitrary (they just
have to satisfy the condition $\Delta<10^{14}$) without the
integration paths getting unnecessarily close to other branch points, which would
affect numerical accuracy. In addition, the integration paths will always be the same
no matter what the relative position of the branch points is. Thus, we consider the
cycles in Fig.~\ref{fig:cutsystem} in the limit that the cycle collapses to a
contour touching the encircled respective cut on both sides,
\begin{equation}
  \label{eq:4}
  \lim_{\epsilon\to0}\int_{C(\epsilon)} \nu_k,
\end{equation}
where $C(\epsilon)$ is any contour which circles around two successive branch
points, such as one of $a_j$ or $c_j$ and where $\epsilon$ is any measure of the
transverse width of the contour across the line between the two branch points. We
assume that in the limit $\epsilon\to0$, this contour degenerates into two lines
between the two branch points. Let us denote by $\Gamma_i^\pm$ the path along the
line between $z_i$ and $z_{i+1}$ lying in the $\pm$-sheet of the Riemann surface,
oriented from $z_i$ to $z_{i+1}$. Assuming the $a$-cycles to lie in the $+$-sheet,
we have
\begin{equation}
  \label{eq:5}
  \lim_{\epsilon\to0} a_j(\epsilon) = \Gamma^+_{2j+1} - \Gamma^-_{2j+1}, \quad \text{and} \quad
  \lim_{\epsilon\to0} c_j(\epsilon) = \Gamma^-_{2j} - \Gamma^+_{2j}.
\end{equation}
Therefore, the contour integrals become 
\[
\int_{a_j}\nu_k = 2\int_{\Gamma_{2j+1}^+}\mskip-20mu\nu_k, \qquad
\int_{c_j}\nu_k = - 2\int_{\Gamma_{2j}^+}\mskip-20mu\nu_k.
\]
Defining 
\begin{equation}
I_{jk} = 2\int_{z_j}^{z_{j+1}}\mskip-20mu\nu_k = 2\int_{\Gamma_j^+}\mskip-10mu\nu_k
\label{eq:7}
\end{equation}
we have ($I_{1k}$ is the first cut, which does not contribute)
\begin{equation}
\label{eq:6}
A_{lk} = I_{2l+1,k}, \qquad B_{lk} = \sum_{j=1}^{l} I_{2j,k}.
\end{equation}

So the problem to compute the periods of the holomorphic differentials $\nu_{k}$,
$k=1,\ldots,g$ is reduced to the computation of some line integrals along a contour
in the complex plane.  Since the $\nu_{k}$ are holomorphic one-forms on the whole Riemann
surface, they are in particular so on the contour. This makes a spectral approach
attractive since it is well known that such approximations have exponential
convergence for analytic functions.

We use here the Clenshaw-Curtis algorithm which is equivalent to an expansion of the
integrand in terms of Chebyshev polynomials. As in the previous subsection, the
contours are mapped to the interval $[-1,1]$ and a set of collocation points
$s_{j}$, $j=0,\ldots,N_{c}$ is introduced on this interval. For Clenshaw-Curtis, the
\emph{Chebyshev points} $s_{j}=\cos (j\pi/N_{c})$, $j=0,\ldots,N_{c}$ are used. To
compute an integral of the form $\int_{-1}^{1}f(s)\mathrm{d}s$, the function $f(s)$
is approximated by Chebyshev polynomials $T_{n}(s)$, $f(s) \approx
\sum_{n=0}^{N_{c}}a_{n}T_{n}(s)$, where the coefficients $a_{n}$ are determined via
a \emph{collocation method}, i.e., by imposing equality of the previous relation at
the collocation points $s_{j}$, $f(s_{j})=\sum_{n=0}^{N_{c}}a_{n}T_{n}(s_{j})$,
$j=0,\ldots,N_{c}$.  Consequently, the integral of $f(s)$ is approximated by
$\int_{-1}^{1}f(s)\,\mathrm{d}s\approx \sum_{n=0}^{N_{c}} a_{n} \int_{-1}^{1}
T_{n}(s)\, \mathrm{d}s = \sum_{n=0}^{N_{c}}w_{n} f(s_{n})$, where the $w_{n}$ are
some known weights (depending on $N_c$) for the Clenshaw-Curtis method (see
\cite{trefethen} and \cite{trefethenweb} for a \Matlab code to compute the
weights). Thus, the integration method consists of sampling the integrand on the
Chebyshev collocation points and computing the scalar product of the vector of
sampled values with the vector built from the weights $w_{n}$. As already mentioned,
this method is known to show exponential convergence for analytic integrands, such
as in the present case of the integration of holomorphic differentials on a
hyperelliptic surface.

The periods are integrals of the form
\[
I = \int_{z_1}^{z_2}\frac{z^k}{\mu(z)}\,\mathrm{d} z
\]
where $\mu(z)$ is the analytically continued square root defined in 
the previous subsection. We compute this integral along the line between $z_1$ and $z_2$ which we parametrize by
\[
z: [-1,1] \to \mathbb{C}, s\mapsto z(s) = v + u s, \qquad 
v:=\frac{z_1+z_2}2, \quad u:=\frac{z_2 - z_1}2. 
\]
Then $\mu(z)$ becomes
\[
\mu(z(s)) = \sqrt{-u^2(1-s^2)p(z(s))}, 
\]
where $p(z) = \prod_{k\ne1,2}(z-z_k)$. Thus the integral reads
\[
I = u\int_{-1}^1 \frac{z(s)^k}{\sqrt{-u^2(1-s^2)p(z(s))}}\,\mathrm{d} 
s
\]
To regularize the integral we substitute $s=\sin(\frac\pi2 t)$ with 
$\mathrm{d} s = \frac\pi2\cos(\frac\pi2 t)$, and the integral becomes
\begin{equation}
 I = \frac\pi2 \int_{-1}^1 \frac{u\cos(\frac\pi2 
t)}{\sqrt{-u^2\cos^2(\frac\pi2 t)p(z(s(t)))}}\,z(s(t))^k\,\mathrm{d} 
t=\pm \mathrm{i}\frac\pi2 \int_{-1}^1 \frac{1}{\sqrt{p(z(s(t)))}}\,z(s(t))^k\,\mathrm{d} 
t.
       \label{sign}
\end{equation}
Since the decision about which sign to use will be made later we take here the positive sign.
To evaluate this integral numerically we use the Clenshaw-Curtis algorithm as
detailed above. The same procedure is applied for all integrals in (\ref{eq:7})
which gives the $(2g+1)\times g$-matrix $I_{ik}$ of line integrals from which we
can form the period matrices.

If one of the integration limits above is infinite, i.e., when $\infty$ is one of
the branch points, then the procedure will be slightly altered. We apply a M\"obius transformation to bring the infinite branch point to a finite location and then compute the appropriately transformed line integrals. In particular, we choose the inversion which interchanges $\infty$ with the central point $c$: $z \mapsto w = 1/(z-c)$. This works except for very symmetric cases when $c$ could be one of the branch points. Then we choose a random point for the inversion. The holomorphic differentials $\nu_k$ are transformed into
\[
\hat\nu_k = -\frac{(1+w c)^{k-1}w^{g-k}}{\tilde\mu(w)}\,\mathrm{d} w, \qquad k=1,\ldots,g
\]
where $\tilde\mu(w)^2 = w^{2g+2} \mu(1/w + c)^2$. After this transformation all branch points lie in a finite place and we compute the line integrals of the holomorphic differentials between infinity and the first branch point by computing the integrals of the transformed differentials $\hat\nu_k$ along the path between $c$ and the image of the first branch point. The other line integrals are computed as described above. 
Alternatively one could compute all the integrals with respect to the 
transformed variables, an approach which we have not explored so far.

Another possibility to address a branching of the surface at infinity 
is to use adapted local coordinates near infinity in the integration, 
i.e.,  $1/\sqrt{z}$ as a local coordinate in an open
neighborhood of infinity. To avoid numerical problems for $z\approx 0$, the
integration path $[z_{j},\infty]$ is split into two intervals $[z_{j},R]$ and
$[R,\infty]$ where $R$ is a point on the chosen ray towards infinity with $|R| \geq
1$ (we take $R$ to be the maximum of 1 and $2|z_{2}|$, the branch point of the cut
extending to infinity). Thus an integral $J$ is computed as a sum of two integrals
$J_{1}$ and $J_{2}$ on these two intervals.  On the first interval, we use
$s=\sqrt{z-z_{j}}$ as a local coordinate, and on the second we take
$\tilde{s}=1/\sqrt{z}$. In both cases we obtain an analytic integrand and use the
Clenshaw-Curtis method as above to compute the integrals. The relative sign of these
two integrals (recall that there is a sign ambiguity due to the square roots) is
fixed by an auxiliary computation to a point $\tilde{z}$ on the integration path of
the original integral with $\tilde{z}\gg R$. The result will not give an accurate
approximation to the integral $J$, but will allow us to fix the sign in the sum
$J_{1}\pm J_{2}$ uniquely. Both of the above approaches are 
implemented in the code and produce identical results within 
numerical accuracy. 

\subsection{Almost degenerate surfaces}

An interesting limit of Riemann surfaces, for instance in the context of
algebro-geometric solutions to integrable equations, is when the surface partially
degenerates and the genus changes by one after pinching a cycle. In the case of a
hyperelliptic surface, this means that two or more branch points coincide in the
limit. In the case of a collapse of a branch cut, it is numerically convenient if
the resulting double point on the curve is surrounded by an $a$-cycle 
(this facilitates the regularization procedure to obtain a Riemann 
surface of genus $g-1$). Since the algorithm to
arrange the branch points in subsection \ref{sec:align-branch-points} does not take
care of almost degenerate situations, as in \cite{lmp} a third option to call the
code exists where the branch points are prearranged to pairs representing the cuts
encircled by the $a$-cycles by the user.

Since the cut-system in Fig.~\ref{fig:cutsystem} is adapted to this case, the
$a$-periods can be treated as before. For the $b$-periods, the fact that the branch
points on the cut crossed by the $b$-cycle are close will affect numerical accuracy.
For illustration, we consider the case that $z_{i}\approx z_{i+1}$. The corresponding
$a$-period is computed as described in the previous subsection, and the integrand in
(\ref{sign}) stays regular even in the limit $z_{i}\to z_{i+1}$. For the integral
between $z_{i+1}$ and $z_{i+2}$ which is needed to compute the $b$-periods much
higher numerical resolution would be needed since $z_{i}\approx z_{i+1}$. To address
this problem, the integral from $z_{i+1}$ to $z_{i+2}$ is split into two integrals
from $z_{i+1}$ to $(z_{i+1}+z_{i+2})/2$ and from $(z_{i+1}+z_{i+2})/2$ to
$z_{i+2}$. In the former case we use for the integration the local coordinate
\[
z = \frac{z_{i}+z_{i+1}}{2}+\frac{z_{i+1}-z_{i}}{2}\cosh s,\quad 
s\in\left[0,\mathrm{arcosh}\left(\frac{z_{i+2}-z_{i}}{z_{i+1}-z_{i}}\right)\right]
\]
or 
$$ s = \sqrt{z-z_{i}}$$ if $z_{i}=\infty$ 
for the former integrand and
\[
z = \frac{z_{i+3}+z_{i+2}}{2}+\frac{z_{i+2}-z_{i+3}}{2}\cosh \tilde{s},\quad 
\tilde{s}\in\left[0,\mathrm{arcosh}\left(\frac{z_{i+1}-z_{i+3}}{z_{i+2}-z_{i+3}}\right)\right]
\]
for the latter. In these local coordinates the integrands can be as well resolved
numerically as the integrands of the $a$-periods even in situations close to the
solitonic limit. After a linear transformation, the integrals are computed with the
Clenshaw-Curtis routine.  Again, the relative sign of these integrals has to be fixed
because of the sign ambiguity in (\ref{sign}). This is done once more by directly
computing the integral between $z_{i+1}$ and $z_{i+2}$ as in the case of the
$a$-periods. This will not be as accurate as the result for the splitting of the
integration path, but it will be enough to fix the relative sign.

Thus, the integrals will be computed efficiently with spectral accuracy, even if two branch
points almost collapse. The sign ambiguity is the same as before and will be
addressed in the following subsection. If more than two branch points almost
collapse as for instance in the so-called positonic limit of KdV, see \cite{pos} and
references therein, the above procedure can only partially address the resulting
numerical problems. As was shown in \cite{pos}, it nonetheless still permits accurate
numerical studies of the situation.

\subsection{Fixing signs via Riemann's bilinear relations}
\label{sec:sorting-out-periods}

Clearly, there is a sign ambiguity in the expression (\ref{sign}) on the right hand
side. This can in principle be fixed by analytically continuing the square root as
described in subsection~\ref{sec:cont-square-root} along a contour close to the
contour given by the cuts in Fig.~\ref{fig:genericcuts}. But this is a
numerically unstable procedure. Instead we always fix the sign only locally along any given cut
(essentially accepting the choice made by \Matlab), and we sort out the signs at
the respective cuts \emph{a posteriori} from the computed period matrices essentially using
Riemann's bilinear identities. Since the cycle $a_{0}$ in Fig.~\ref{fig:cutsystem}
is homologous to a sum of the cycles $a_{j}$, $j=1,\ldots,g$, this implies a
condition on the periods $A_{kl}$ in (\ref{eq:6}).  This is a $(g+1)\times
g$-matrix, which is the matrix of $a$-periods when we discard the first
row. However, each column of this matrix must add up to zero if the root $\mu(z)$ is
analytically continued according to Fig.~\ref{fig:cutsystem}. Since this is not enforced
in our construction, we only know that there is a left null vector $n^i$, so that
\[
n^iA_{ik} = 0.
\]
If the function $\mu$ was analytically continued then this vector will be proportional to $(1,1,\ldots,1)$, while there will be negative signs otherwise. In order to find the correct signs of the $n^i$, we identify the null space of
$A_{ik}$. \Matlab does this  by performing a
\emph{singular value decomposition} (SVD) which is for an
  $m\times n$-matrix $M$ with complex entries given by $M=\mathcal{U}\Sigma
  \mathcal{V}^{\dagger}$; here $\mathcal{U}$ is an $m\times m$ unitary matrix,
  $\mathcal{V}^{\dagger}$ denotes the conjugate transpose of $\mathcal{V}$, an $n\times n$
  unitary matrix, and the $m\times n$ matrix $\Sigma$ is diagonal (as
  defined for a rectangular matrix); the non-negative numbers on the
  diagonal of $\Sigma$ are called the \textit{singular values} of
  $M$. The line in $\mathcal{U}$ corresponding to the singular value with 
  smallest modulus gives the vector $n^{i}$. If sufficient numerical resolution is
provided, i.e., if $N_{c}$ is large enough, this null space is 1-dimensional and can
be normalized to be a vector with elements $\pm 1$, i.e., $n^i = \pm1$. A failure of this
procedure implies that $N_{c}$ was too small, and that the code has to be rerun with
a larger value of $N_{c}$ to compute the periods with higher precision. Again, the
spectral convergence of the code is very useful here. Once the vector with
components $n^{i}$ is correctly identified, it consists of $\pm1$'s and we now
replace $A_{ik} \leftarrow n^i A_{ik}$ (no summation over $i$) and hence fix the
signs of the $a$-periods. The first row of $A_{jk}$ can then be 
discarded, and the resulting matrix is denoted by~$\mathcal{A}$. 

In order to fix the signs of the $b$-periods we proceed in a similar way. Note, that
the real part of the Riemann matrix
\begin{equation}
    \mathbb{B} = 2\pi \mathrm{i}\, \mathcal{B} \mathcal{A}^{-1}
   \label{RieMat},
\end{equation}
where $\mathcal{B}$ is the matrix of $b$-periods of the holomorphic 
one-forms, must be symmetric and negative definite. The signs of those line
integrals which make up the $a$-periods are already correct. Now we construct the
$b$-periods according to (\ref{eq:6}).   We allow for possibly different signs by
writing
\[
B_{lk} = - \sum_{i=1}^l \epsilon_i I_{2i,k}
\]
for $g$ undetermined signs $\epsilon_i$. Computing the matrix $\mathcal{B}$ and then the
Riemann matrix with this expression, we can write down the symmetry condition on
$\mathbb{B}$
\[
\operatorname{Re}(\mathbb{B} - \mathbb{B}^t) = 0.
\]
These are $\binom{g}{2}$ conditions, which are linear in the $\epsilon_i$.  We choose the
$g-1$ conditions which come from the first off-diagonal of this matrix equation. This
yields $g-1$ linear equations for the $g$ unknowns $\epsilon_i$. Thus, generically, we will
again get a system with a 1-dimensional kernel which is identified as explained
for the $a$-periods. This leads again to a null vector which has as entries only $\pm1$
(up to numerical error of course) which give the relative signs of the line
integrals to make the upper and lower off-diagonals in the real part of the Riemann
matrix equal. We fix the remaining overall sign by the requirement 
that the real part of $\mathbb{B}$ should be negative definite. 
This completely fixes the signs of the
line integrals. 

\begin{figure}[htb]
  \centering
  \includegraphics[width=0.7\linewidth]{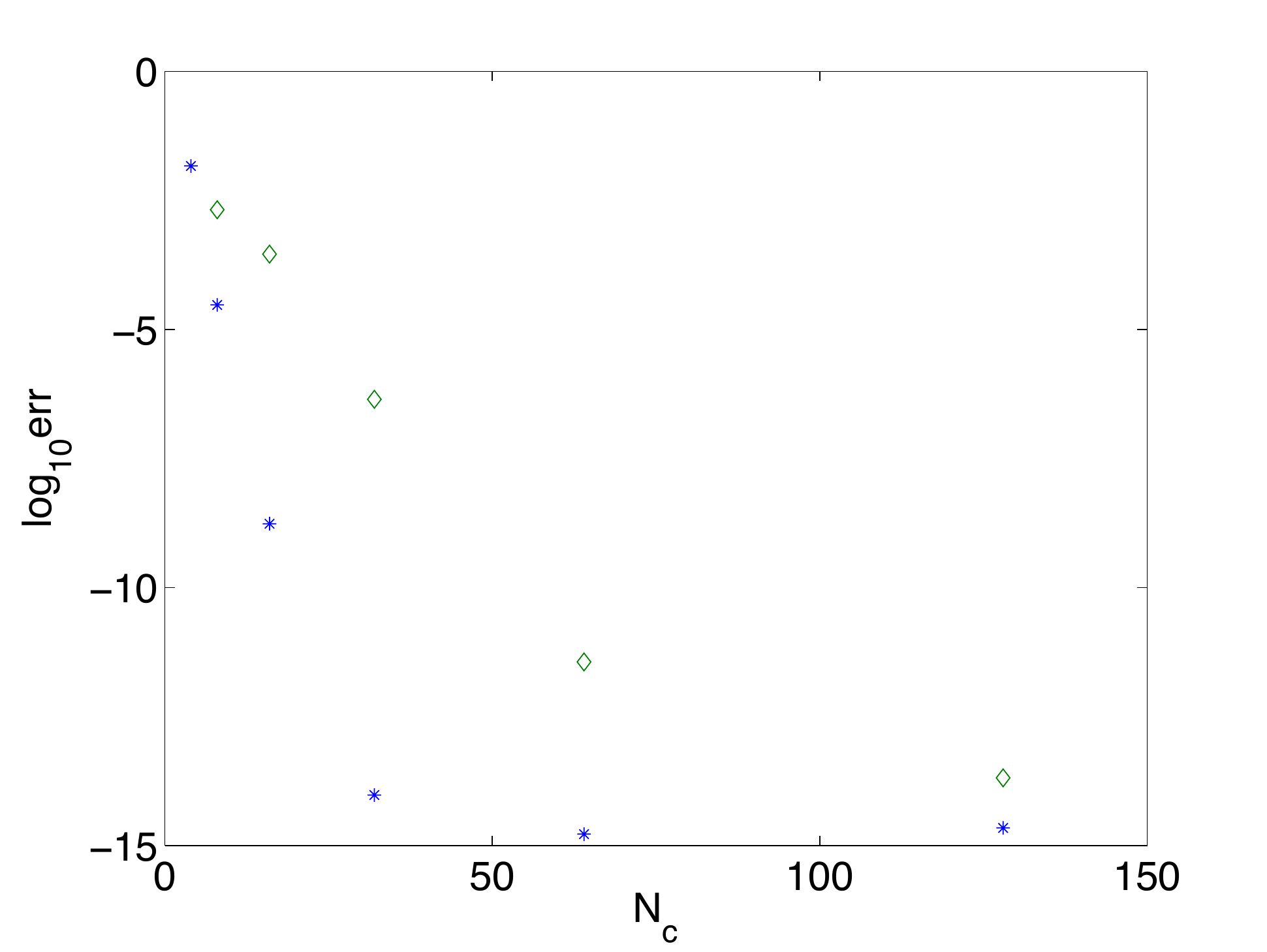}
  \caption{Numerical error $err$ defined as the maximum of the 
  $L^{\infty}$ norm of the skew 
  part of the Riemann matrix and the $L^{\infty}$ norm of the sum of 
  the $a$ periods including the period over $a_{0}$ in dependence of 
  the number of Chebyshev points $N_{c}$ for the example of the genus~2 curve given by
  $\mu^{2}=z(z-1-\mathrm{i}\epsilon)(z-1+\mathrm{i}\epsilon)(z-2-\mathrm{i}\epsilon)(z-2+\mathrm{i}\epsilon)$; 
  the stars correspond to the case $\epsilon=1$, the diamonds to $\epsilon=10^{-12}$.}
  \label{hypererr}
\end{figure}

The size of the remaining skew part of the real part of the Riemann
matrix provides a good test of the accuracy of the numerical 
approach. We use the maximum of the $L^{\infty}$ norm of this skew 
part and the difference between the sum of the $a$-periods and the 
periods of the holomorphic one-forms along the cycle $a_{0}$ in 
Fig.~\ref{fig:genericcuts}, denoted by $err$ in the following,  as  
an indicator of the accuracy of the computation. A typical example 
for this can be seen in Fig.~\ref{hypererr}. Visibly the 
error decreases exponentially with the number $N_{c}$ of Chebyshev 
points, i.e., we get the expected spectral convergence. This also 
holds in an almost degenerate situation as can be seen in the same 
figure, where pairs of branch points are separated only by $10^{-12}$ 
in the second example. In both cases machine precision is reached 
with $N_{c}=128$. It is remarkable that in the first case, $N_{c}=4 
$ already leads to an error of the order of a few percent. As 
mentioned, the code is very efficient. On an average computer, the 
computation of the periods takes of the order of 1ms in the shown 
example for $N_{c}=128$. Note that the maximally achievable accuracy 
for the Riemann matrix in general depends on the conditioning (the 
range of the eigenvalues) of the matrix $\mathcal{A}$ of $a$-periods 
since the latter has to be inverted in (\ref{RieMat}) to determine 
$\mathbb{B}$. This limits the achievable accuracy for the Riemann matrix for higher 
genus even if all periods are computed to the order of machine 
precision.

\subsection{Abel map}

Integrals of the holomorphic one-forms between arbitrary points $P$, $P_{0}$ of
the hyperelliptic Riemann surface $\mathcal{R}$ can be computed essentially in the same way as
the periods of the holomorphic differentials above.  The \emph{Abel map}
\[
P \mapsto \int_{P_{0}}^{P}\nu_k 
\]
is a bijective map from the surface into the
\emph{Jacobian} $Jac(\mathcal{R}) = \mathbb{C}^{g}/\Lambda$, where $\Lambda$ is
the lattice formed by the periods of the holomorphic one-forms, 
\[
\Lambda =
\left\{2\pi \mathrm{i}\,\mathrm{m} +\mathbb{B} \,\mathrm{n}: \mathrm{m}, \mathrm{n}\in \mathbb{Z}^{g}\right\}.
\]

To compute this map for a given point $P\in \mathcal{R}$, we identify the branch
point closest to $P$, which will be denoted by $E$. Then we compute the integrals
$\int_{E}^{P}\nu_k$ (we only discuss here the computation of the integral
between $E$ and $P$, the one between $P_{0}$ and $E$ can then be computed in the
same way). It is a consequence of the relations (\ref{eq:6}) that the Abel map
between branch points on a hyperelliptic surface are half-periods. Thus, if $P$
is a branch point, the Abel map has been already computed above. If $P$ is a
finite point, we use $s=\sqrt{z-E}$ as a local coordinate and compute the line
integral after the analytic continuation of the square root as in
\ref{sec:cont-square-root} with the Clenshaw-Curtis algorithm applied to the
integral in $s$.

If $P$ is a point covering infinity on a surface without branching at infinity,
we introduce some intermediate point $Q$ with $\overline{QE}\gg 1$. The integral
to $Q$ is computed as described for a finite $P$. For the integral between $Q$ and
$P$, we use $s=1/z$ as a local parameter and compute the resulting integral
again with the Clenshaw-Curtis method.

Since the Abel map is only defined up to periods of the holomorphic one-forms, we
always choose it to be in the fundamental domain given by
\[
 2\pi \mathrm{i}\,\mathrm{p} + \mathbb{B}\,\mathrm{q},\quad \mathrm{p},\mathrm{q}\in \mathbb{R}^{g},
\]
with $-1/2<p_{i}\leq 1/2$, $-1/2<q_{i}\leq 1/2$, $i=1,\ldots,g$. In 
other words, the Abel map and thus an arbitrary point of the Jacobian 
can be given in terms of the $(p_{i}, q_{i})$, which are called 
\emph{characteristics}. Of special importance are half-integer 
characteristics with $2\mathrm{p},2\mathrm{q}\in \mathbb{Z}^{g}$.  A 
half-integer characteristic is called \emph{even} if $4\langle \mathrm{p},\mathrm{q}\rangle=0\mbox{ mod } 2$ and \emph{odd} otherwise. Here  
    $\left
    \langle\cdot,\cdot\right\rangle$ denotes the Euclidean scalar product
    $\left\langle \mathrm{p},\mathrm{q}\right\rangle=\sum_{i=1}^g 
    p_iq_i$. The Abel map between branch points corresponds thus to a half integer characteristic.

\section{Sine-Gordon equation}
In this section we discuss solutions to the SG equation (\ref{SG}) in 
terms of multi-dimensional theta functions on hyperelliptic Riemann 
surfaces which were first constructed by Kozel and Kotlyarov \cite{kotlyarov}, 
see also \cite{algebro} and \cite{tata} for an alternative 
derivation based on Fay's trisecant identity. We present examples for 
real, smooth solutions, also in 
the solitonic limit. The solutions are given on Riemann surfaces with 
cuts between real points and/or cuts between complex conjugate 
points. Contrary to the previous code, the new version can handle 
surfaces with both types of cuts on the  same surface which will be illustrated in this 
section. Note that the SG equation can be written in the form of a $1+1$ non-linear wave equation
\begin{equation}
    \phi_{tt}-\phi_{xx}=\sin(\phi),\quad \xi=\frac{x+t}{4}, 
    \eta=\frac{t-x}{4}.
    \label{SG2}
\end{equation} 
We will discuss the solutions here only in dependence of the 
characteristic coordinates $\xi$ and $\eta$. 

\subsection{Solutions in terms of multi-dimensional theta functions}

Solutions to the SG equation are given in terms of multi-dimensional 
theta functions which we define 
as an infinite series,
    \begin{equation}\label{theta}
    \Theta_{\mathrm{pq}}(\mathrm{z},\mathbb{B})=
    \sum\limits_{\mathrm{N}\in\mathbb{Z}^g}\exp\left\{\frac{1}{2}
    \left\langle\mathbb{B}\left(\mathrm{N+p}\right),
    \mathrm{N+p}
    \right\rangle+
    \left\langle \mathrm{z+q},\mathrm{N+p}
    \right\rangle\right\}
    \;,
    \end{equation}
    with $\mathrm{z}\in\mathbb{C}^g$ and $\mathrm{p}$, $\mathrm{q}\in{
    \mathbb{R}}^g$ the characteristics.

The properties of the Riemann matrix ensure that the series converges
absolutely and that the theta function is an entire function on
$\mathbb{C}^{g}$.  A characteristics is called \emph{singular} if the
corresponding theta function vanishes identically.  Theta functions 
with odd (even) characteristics are odd
(even) functions of the argument $\mathrm{z}$.  The theta function with
characteristics is related to the Riemann theta function $\Theta$, the
theta function with zero characteristics $\Theta:= \Theta_{00}$,
via
\begin{equation}
    \Theta_{\mathrm{pq}}(\mathrm{z},\mathbb{B})=\Theta(\mathrm{z}
    +\mathbb{B}\mathrm{p} + \mathrm{q})\exp\left\{\frac{1}{2}
    \left\langle\mathbb{B}\mathrm{p,p}\right\rangle+
    \left\langle \mathrm{p,z} + \mathrm{q}\right\rangle
    \right\}\;.
    \label{thchar}
\end{equation}
The theta function has the periodicity properties 
\begin{equation}
    \Theta_{\mathrm{pq}}(\mathrm{z}+2\pi \mathrm{i} e_{j}) = 
    e^{2\pi ip_{j}}
    \Theta_{\mathrm{pq}}(\mathrm{z})\;,
    \quad 
    \Theta_{\mathrm{pq}}(\mathrm{z}+\mathbb{B}
    e_{j})=
    e^{- (z_{j}+q_{j}) - \frac{1}{2} B_{jj}}
    \Theta_{\mathrm{pq}}(\mathrm{z})\;
    \label{eq:periodicity},
\end{equation}
where $e_{j}$ is a vector in $\mathbb{R}^{g}$ consisting of
zeros except for a 1 in jth position.  In the computation of the
theta function, the series~(\ref{theta}) is 
approximated as a sum. The argument $\mathrm{z}$ is always written as 
$\mathrm{z}=\mathrm{z}_{0}+2\pi \mathrm{i}\mathrm{m}+\mathbb{B}\mathrm{n}$, where 
$z_{0}$ is in the fundamental domain of the Jacobian. The theta 
function is computed for the argument $\mathrm{z}_{0}$, and the 
relations (\ref{eq:periodicity}) then give the theta function for the 
argument $\mathrm{z}$.  For details of the computation, the reader is 
referred to \cite{lmp} and \cite{RSbookfk}. 

Solutions to the SG equation can be obtained on a hyperelliptic 
Riemann surface  given by
\begin{equation}
    \mu^{2}=z\prod_{i=1}^{2g}(z-z_{i}),\quad \Re z_{i}>0,
    \label{SGsurf}
\end{equation}
in the form 
\begin{equation}
    u = 2\mathrm{i}\ln \frac{\Theta_{\mathrm{pq}}(\mathrm{V}\xi+\mathrm{W}\eta+\mathrm{i}\pi \delta)}{\Theta_{\mathrm{pq}}(V\xi+W\eta)}
    \label{SGformula},
\end{equation}
where $\mathrm{V,W}\in \mathbb{C}^{g}$ with components
$$
V_{j}= 2c_{j1},\quad 
W_{j}=\frac{2c_{jg}}{\sqrt{\sum_{i=1}^{2g}z_{i}}},\quad j=1,\ldots,g,
$$ with $c_{ik}=2\pi \mathrm{i}(\mathcal{A}^{-1})_{ik}$, $i,k=1,\ldots,g$, 
where $
\begin{bmatrix}
    \mathrm{p}\\
    \mathrm{q}
\end{bmatrix}
$ is an arbitrary nonsingular characteristics, and where 
$\delta\in\mathbb{C}^{g}$ with components $\delta_{i}=1$ for 
$i=1,\ldots,g$. 

The solutions are real and regular (without poles), see the 
discussion in \cite{algebro} (note that a different homology basis is 
used there), if the hyperelliptic surface is real, 
i.e., if the branch points $z_{i}$ are real or pairwise complex 
conjugate, and if the characteristic satisfies certain reality 
conditions. Let $\tau$ be the antiholomorphic involution of this 
surface acting on each sheet as the complex conjugation. 
On such a surface it is possible to introduce a basis of 
the homology of the form shown in Fig.~\ref{fig:cutsystem}
satisfying
$
\tau a_{j}=-a_{j}$, $j=1,\ldots,g$ and $\tau b_{j}=b_{j}$ if $b_{j}$ 
crosses a cut of real branch points only, and $\tau 
b_{j}=b_{j}-a_{j}$ if $b_{j}$ passes through a cut between conjugate 
branch points, see \cite{algebro}. Note that we do not always use such a cut system in 
the following since we also want to consider cuts collapsing away 
from the real axis. For the examples we will study below, we will always 
give the used characteristics and the branch points in pairs 
corresponding to the cuts encircled by the $a$-cycles.

Note that the solution (\ref{SGformula}) is only defined up to 
multiples of $2\pi$. Since the logarithm in \Matlab is branched on 
the negative real axis, solutions computed according to formula 
(\ref{SGformula}) will in general have jumps. As for the analytic 
continuation of the square root discussed in the previous section, we 
construct an analytic solution by comparing $u(\xi,\eta)$ on 
neighboring computed points $|u(\xi_{j},\eta)-u(\xi_{j+1},\eta)|$ to 
$|u(\xi_{j},\eta)-u(\xi_{j+1},\eta)\pm 4\pi|$. If the latter 
expression is smaller than the former, a factor of $\pm4\pi$ is 
added. The same procedure is repeated for the coordinate $\eta$ so that we end up with a 
smooth solution $u$.

The numerical accuracy of the solution is controlled in two ways. 
First we check for each computed point $(\xi,\eta)$ that the  
identity for theta functions obtained by entering (\ref{SG}) with the 
solution $u$ of (\ref{SGformula}),
\[
\begin{split}
&\frac{\partial_{VW}\Theta_{\mathrm{pq}}(\mathrm{V}\xi+\mathrm{W}\eta+\mathrm{i}\pi \delta)}
{\Theta_{\mathrm{pq}}(\mathrm{V}\xi+\mathrm{W}\eta+\mathrm{i}\pi \delta)}-
\frac{\partial_{V}\Theta_{\mathrm{pq}}(\mathrm{V}\xi+\mathrm{W}\eta+\mathrm{i}\pi \delta)\partial_{W}\Theta_{\mathrm{pq}}(\mathrm{V}\xi+\mathrm{W}\eta+\mathrm{i}\pi \delta)
}{\Theta_{\mathrm{pq}}(\mathrm{V}\xi+\mathrm{W}\eta+\mathrm{i}\pi 
\delta)^2} \\
&-\frac{\partial_{VW}\Theta_{\mathrm{pq}}(\mathrm{V}\xi+\mathrm{W}\eta)}
{\Theta_{\mathrm{pq}}(\mathrm{V}\xi+\mathrm{W}\eta)}+
\frac{\partial_{V}\Theta_{\mathrm{pq}}(\mathrm{V}\xi+\mathrm{W}\eta)
\partial_{W}\Theta_{\mathrm{pq}}(\mathrm{V}\xi+\mathrm{W}\eta)
}{\Theta_{\mathrm{pq}}(\mathrm{V}\xi+\mathrm{W}\eta)^2}= \frac{1}{2\mathrm{i}}\sin u
\end{split}
\]
where the
derivative of a function $F(\mathrm{z})$ with respect to a vector $U$ 
is defined 
as 
$\partial_{U}F(\mathrm{z}):=    \left\langle\nabla F(\mathrm{z}),U
    \right\rangle
$.
This identity for theta functions 
is not built into the code and thus provides 
a strong test. The code reports a warning if the identity is not 
satisfied to better than $10^{-12}$.

In addition, the solutions $u$ are computed on Chebyshev collocation points for
$\xi$ and $\eta$. As in the
Clenshaw-Curtis integration in the previous section, this can be used to
approximate the computed solution via Chebyshev polynomials. Since the
derivatives of the Chebyshev polynomials can be expressed linearly in terms of
Chebyshev polynomials, a derivative acts on the space of polynomials via a
\emph{differentiation matrix}. With these standard Chebyshev differentiation
matrices (see \cite{trefethen,trefethenweb}), the solution can be numerically
differentiated. With the computed derivatives we check to which numerical
precision the partial differential equation (PDE) is satisfied by the numerical
solution. If sufficient resolution is provided, the residual of the equation
obtained with Chebyshev differentiation is smaller than $10^{-6}$.

\subsection{Examples}

We first study examples on real surfaces of genus 2. The branch 
points are given in pairs corresponding to the chosen cuts. The first 
cut is always between $-\infty$ and 0, the remaining cuts are freely 
chosen according to the reality conditions. 
We always use $N_{c}=128$ Chebyshev points unless 
otherwise noted. 

\begin{figure}[htb] 
  \includegraphics[width=0.49\linewidth]{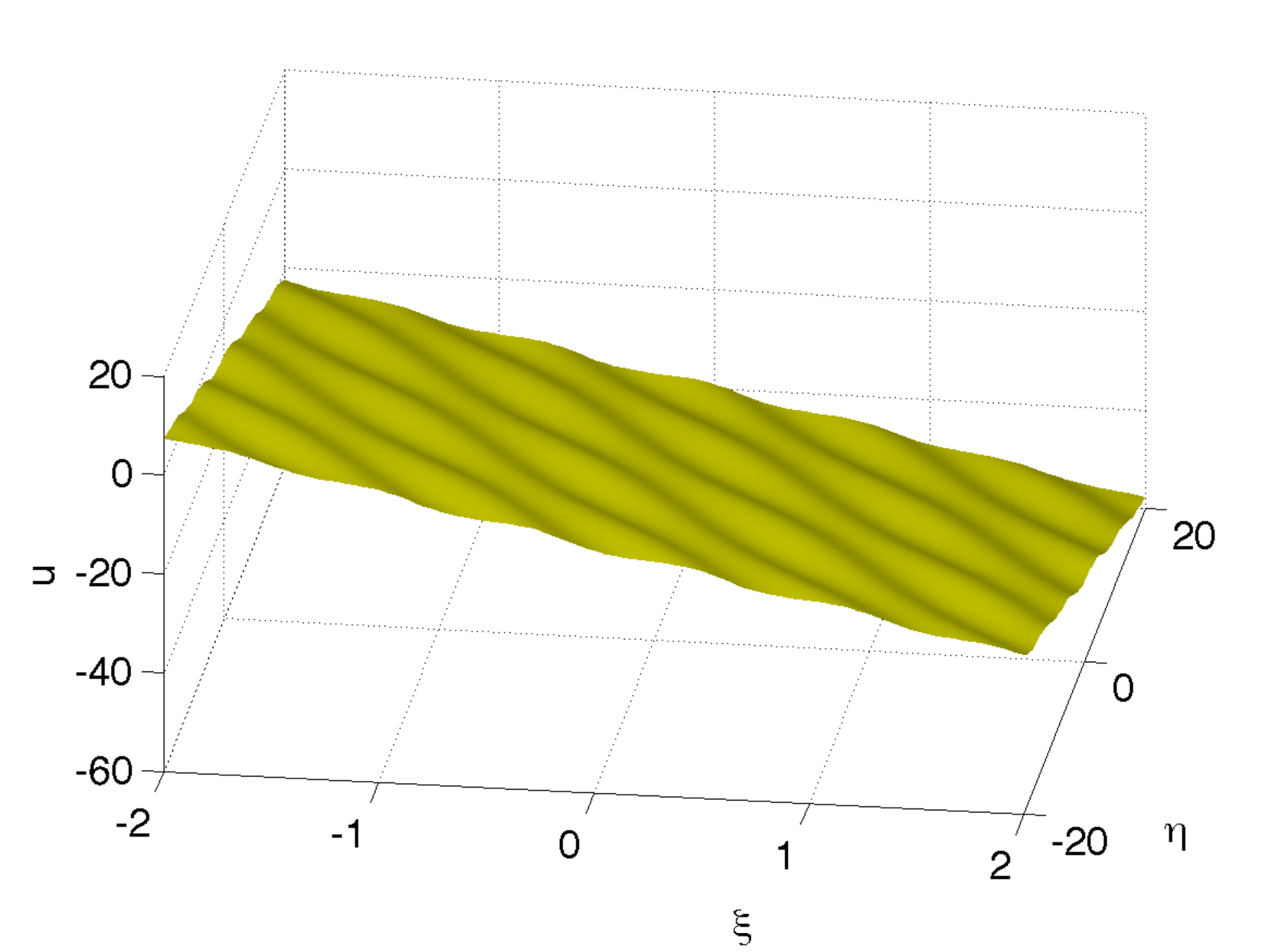}
  \includegraphics[width=0.49\linewidth]{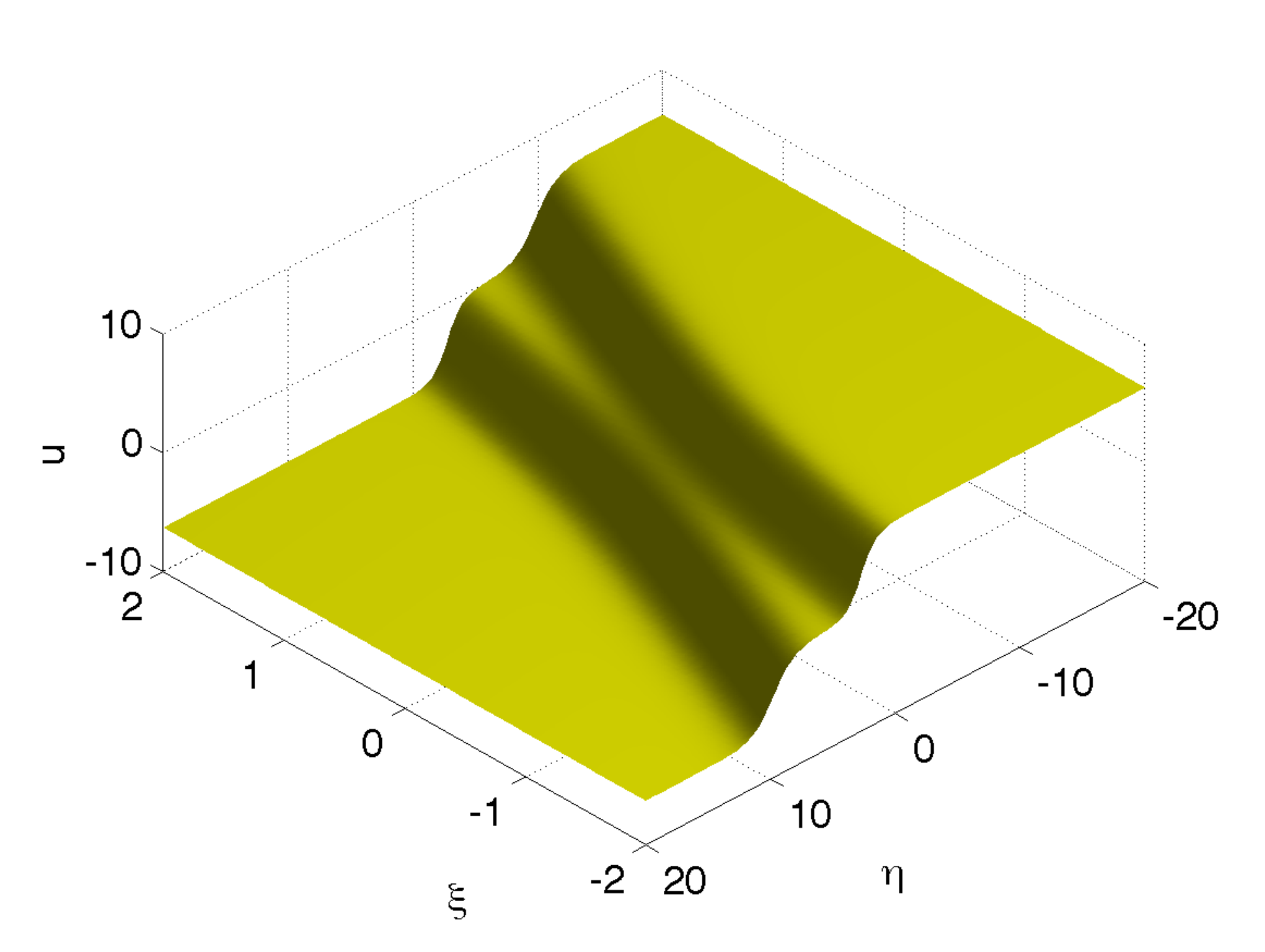}
  \caption{\label{SGrere} Solution (\ref{SGformula}) to the SG equation (\ref{SG}) 
  on the surface of genus 2 with branch points (\ref{rere}), on the left for $\epsilon=1$, on the right for 
  $\epsilon=10^{-12}$. }
\end{figure}

On the 
surface with the real branch points 
\begin{equation}
   \begin{bmatrix}
      -\infty & 1 & 3\\
      0 & 1+\epsilon & 3+\epsilon
  \end{bmatrix}
       \label{rere}
\end{equation}
we use the characteristics $  \begin{bmatrix}
            \mathrm{p} \\
      \mathrm{q}
  \end{bmatrix}
=\frac{1}{2}
  \begin{bmatrix}
            1& 1 \\
      1/2 & 1/2
  \end{bmatrix}
  $.
As can be seen in Fig.~\ref{SGrere}, one obtains for 
$\epsilon=1$ a quasiperiodic solution. For $\epsilon\to0$, this 
solution becomes the 2-kink solution of the SG equation as can be 
seen in the same figure on the right. Thus the solution on 
  the non-degenerate surface can be seen as an infinite train of such 
  2-kinks. 

\begin{figure}[htb]  
  \includegraphics[width=0.49\linewidth]{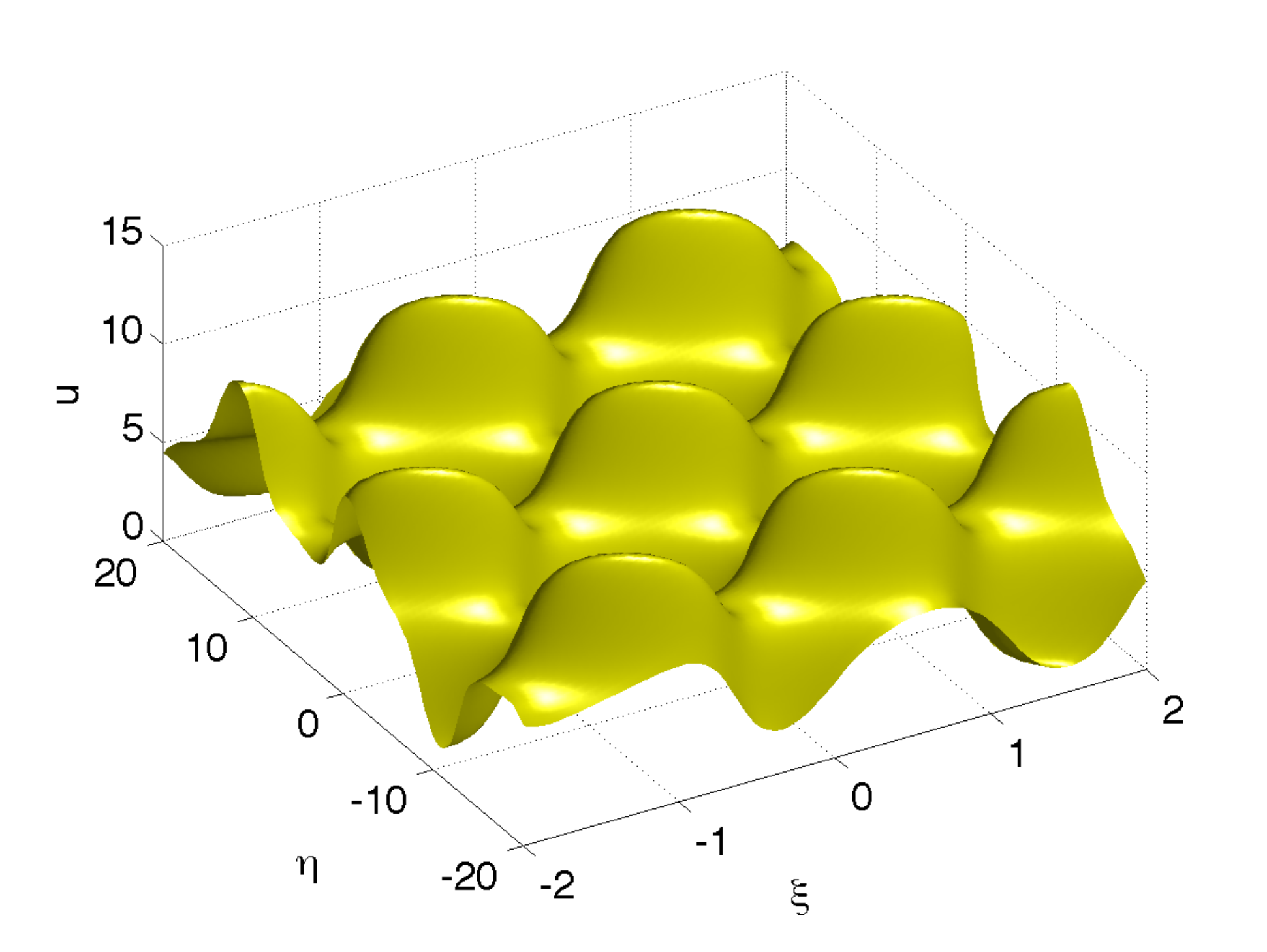}
  \includegraphics[width=0.49\linewidth]{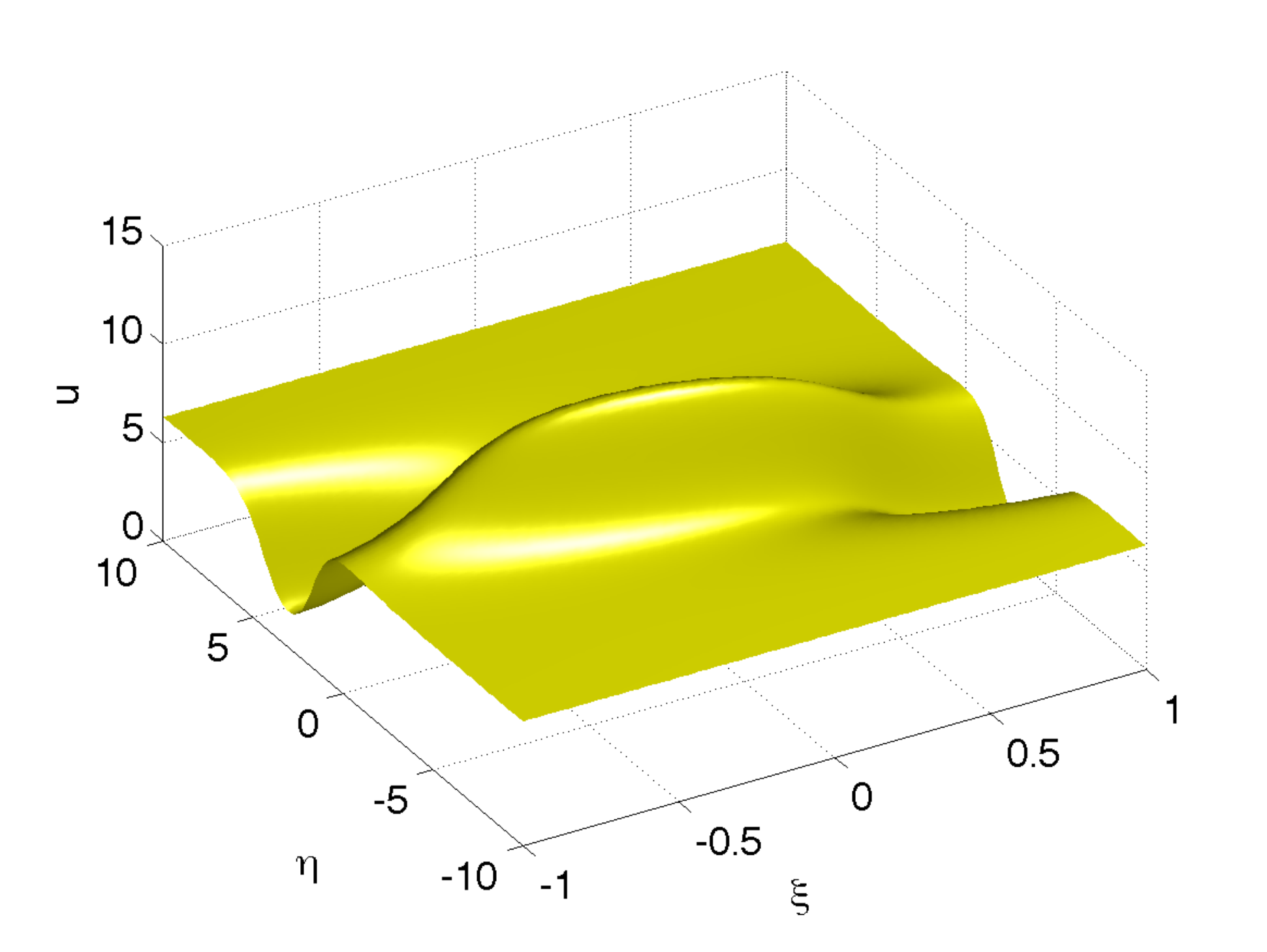}
  \caption{\label{SGim} Solution (\ref{SGformula}) to the SG equation (\ref{SG}) 
  on the surface with branch points (\ref{imim}) ; on the left for $\epsilon=1$, on the right for 
  $\epsilon=10^{-12}$. }
\end{figure}

On the surface with the cuts 
\begin{equation}
    \begin{bmatrix}
      -\infty & 1+\mathrm{i} &  1-\mathrm{i}\\
      0 & 1+\mathrm{i}+\epsilon & 1-\mathrm{i}+\epsilon
  \end{bmatrix}
      \label{imim},
\end{equation}
i.e., only conjugate branch points except for $[-\infty,0]$,
we use the characteristics $  \begin{bmatrix}
            \mathrm{p} \\
      \mathrm{q}
  \end{bmatrix}
=\frac{1}{2}
  \begin{bmatrix}
            1& 1 \\
      0 & 0
  \end{bmatrix}
  $. In the limit $\epsilon\to 0$, double points appear this time not 
  on the real axis, but are conjugate to each other.  The solutions 
  for this case can be seen in Fig.~\ref{SGim} for $\epsilon=1$ on 
  the left  and on 
  an almost degenerate surface ($\epsilon=10^{-12}$) on the right. 

\begin{figure}[htb]  
  \includegraphics[width=0.7\linewidth]{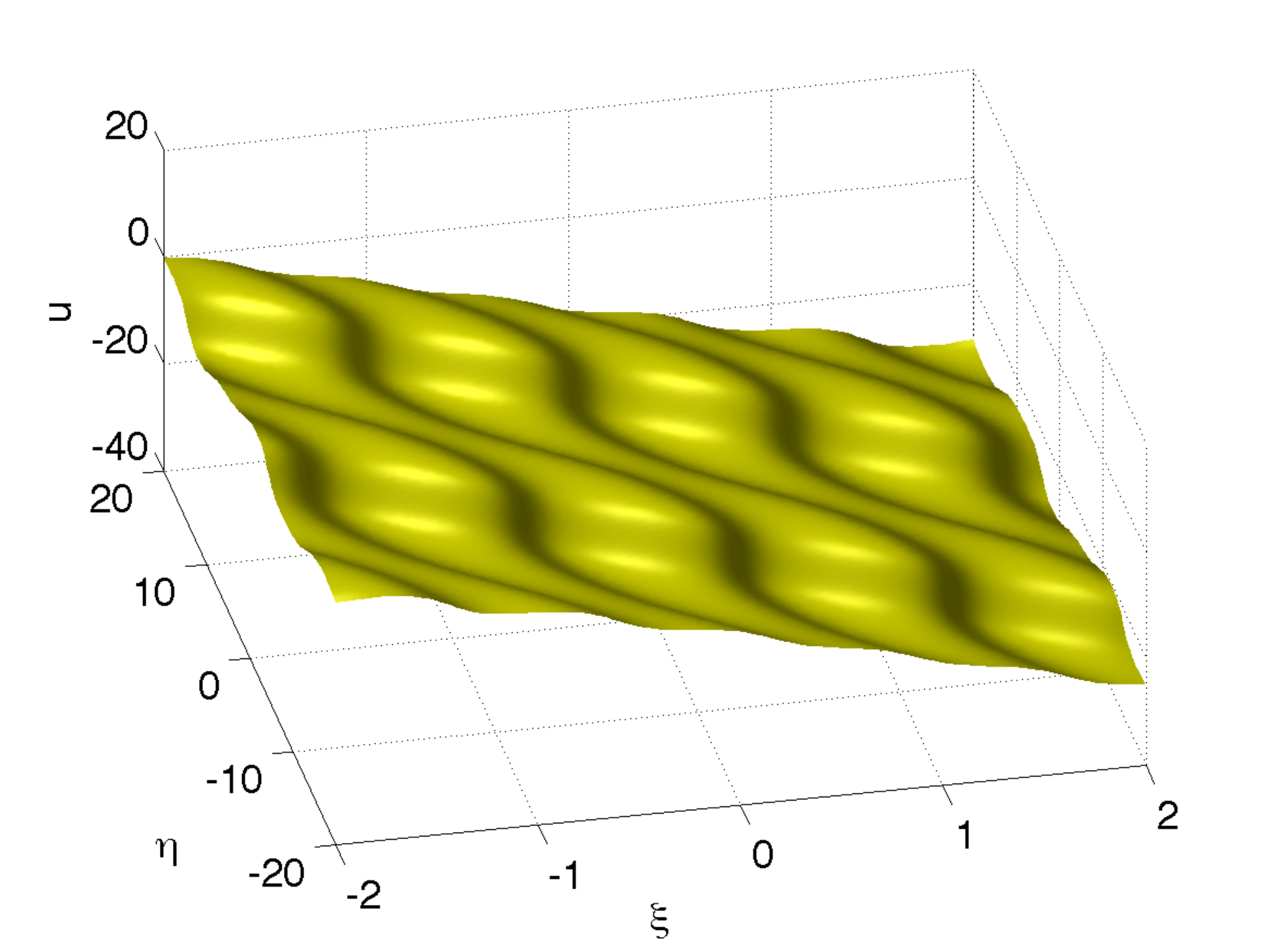}
  \caption{\label{SGreim} Solution (\ref{SGformula}) to the SG equation (\ref{SG}) 
  on the surface with branch points (\ref{reim}). }
\end{figure}

An example for mixed real and conjugate cuts in addition to the cut $[-\infty,0]$
is given on the genus 2 surface with the 
cuts
\begin{equation}
     \begin{bmatrix}
      -\infty & 1 &  3+\mathrm{i}\\
      0 & 2 & 3-\mathrm{i}
  \end{bmatrix}
     \label{reim},
\end{equation}
for which we use the characteristics $  \begin{bmatrix}
           \mathrm{ p} \\
      \mathrm{q}
  \end{bmatrix}
=\frac{1}{2}
  \begin{bmatrix}
            1& 1 \\
      1/2 & 1/2
  \end{bmatrix}
  $. The corresponding solution can be seen in Fig.~\ref{SGreim}.

\begin{figure}[htb]
  \includegraphics[width=0.49\linewidth]{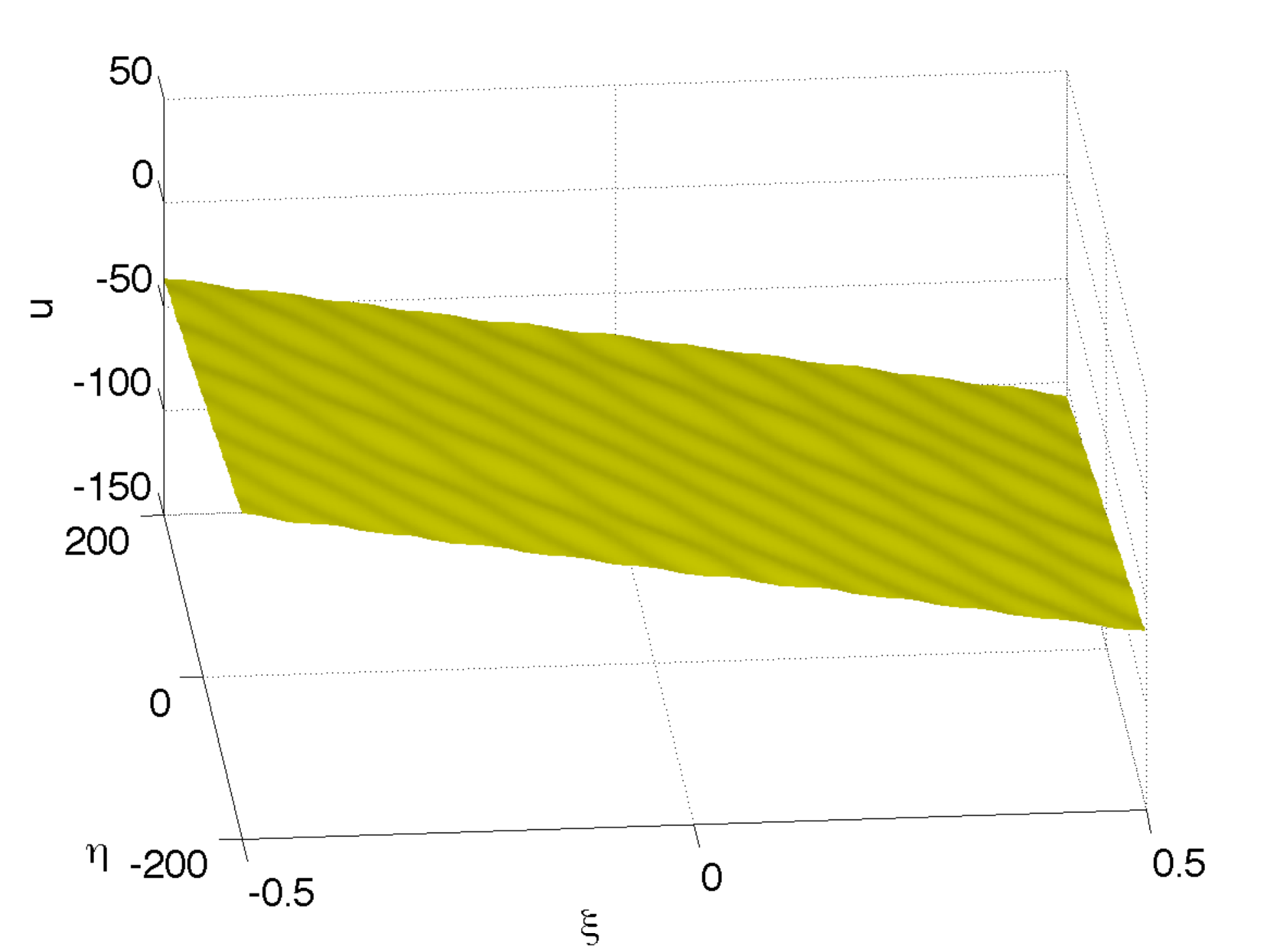}
  \includegraphics[width=0.49\linewidth]{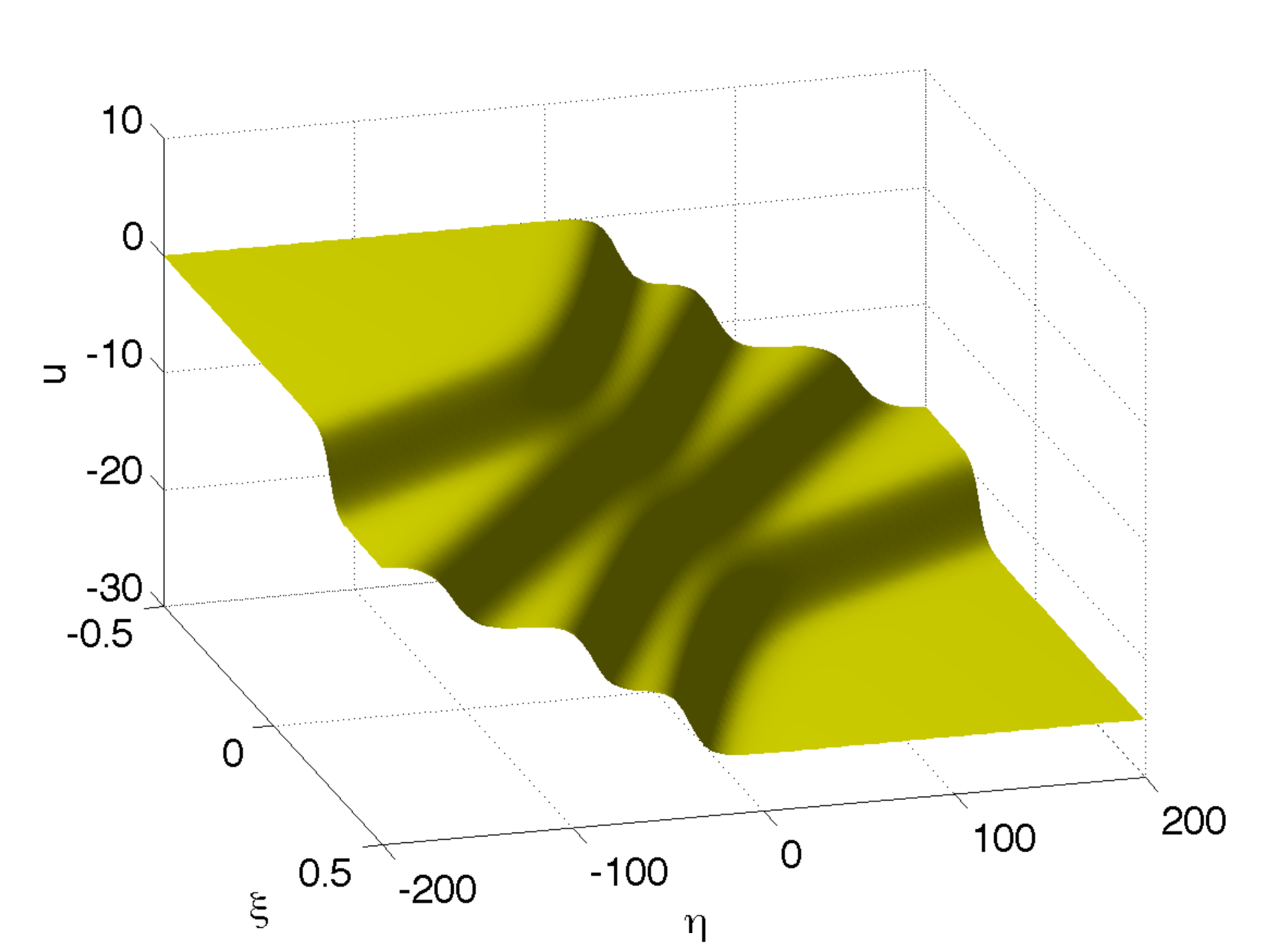}
  \caption{\label{SGreg4} Solution (\ref{SGformula}) to the SG equation (\ref{SG}) 
  on the surface with branch points (\ref{reg4}); on the left for 
  $\epsilon=0.5$, on the right for 
  $\epsilon=10^{-12}$. }
\end{figure}

In Fig.~\ref{SGreg4}, we show an example on a hyperelliptic Riemann 
surface of genus 4 with only real branch points,
\begin{equation}
     \begin{bmatrix}
      -\infty & 1 &  2& 3& 4\\
      0 & 1+\epsilon & 2+\epsilon& 3+\epsilon&4+\epsilon
  \end{bmatrix}
  \label{reg4}.
\end{equation}
We use the characteristics $  \begin{bmatrix}
            \mathrm{p} \\
      \mathrm{q}
  \end{bmatrix}
=\frac{1}{2}
  \begin{bmatrix}
            1& 1 & 1 &1\\
      1/2 & 1/2 & 1/2 & 1/2
  \end{bmatrix}
  $. The 4-kink solution appearing in the solitonic limit on the 
  right can be clearly recognized. The quasiperiodic solution on the 
  non-degenerate surface can be seen as an infinite sequence of such 
  kinks.

A further example for a 
surface of genus 4 has the branch points,
\begin{equation}
     \begin{bmatrix}
      -\infty & 1+\mathrm{i} &  1-\mathrm{i}& 2+\mathrm{i}& 2-\mathrm{i}\\
      0 & 1+\mathrm{i}+\epsilon & 1-\mathrm{i}+\epsilon& 2+\mathrm{i}+\epsilon&2-\mathrm{i}+\epsilon
  \end{bmatrix}
  \label{img4},
\end{equation}
i.e., except for $[-\infty,0]$ only  conjugate cuts forming double points away from the 
real axis in the limit $\epsilon\to0$. 
We use the characteristics $  \begin{bmatrix}
            \mathrm{p} \\
      \mathrm{q}
  \end{bmatrix}
=\frac{1}{2}
  \begin{bmatrix}
            1& 1 & 1 &1\\
      0 & 0&0&0
  \end{bmatrix}
  $. To compute the case with $\epsilon=10^{-12}$, we use $N_{c}=512$ 
  Chebyshev points. 
 \begin{figure}[htb]
  \includegraphics[width=0.49\linewidth]{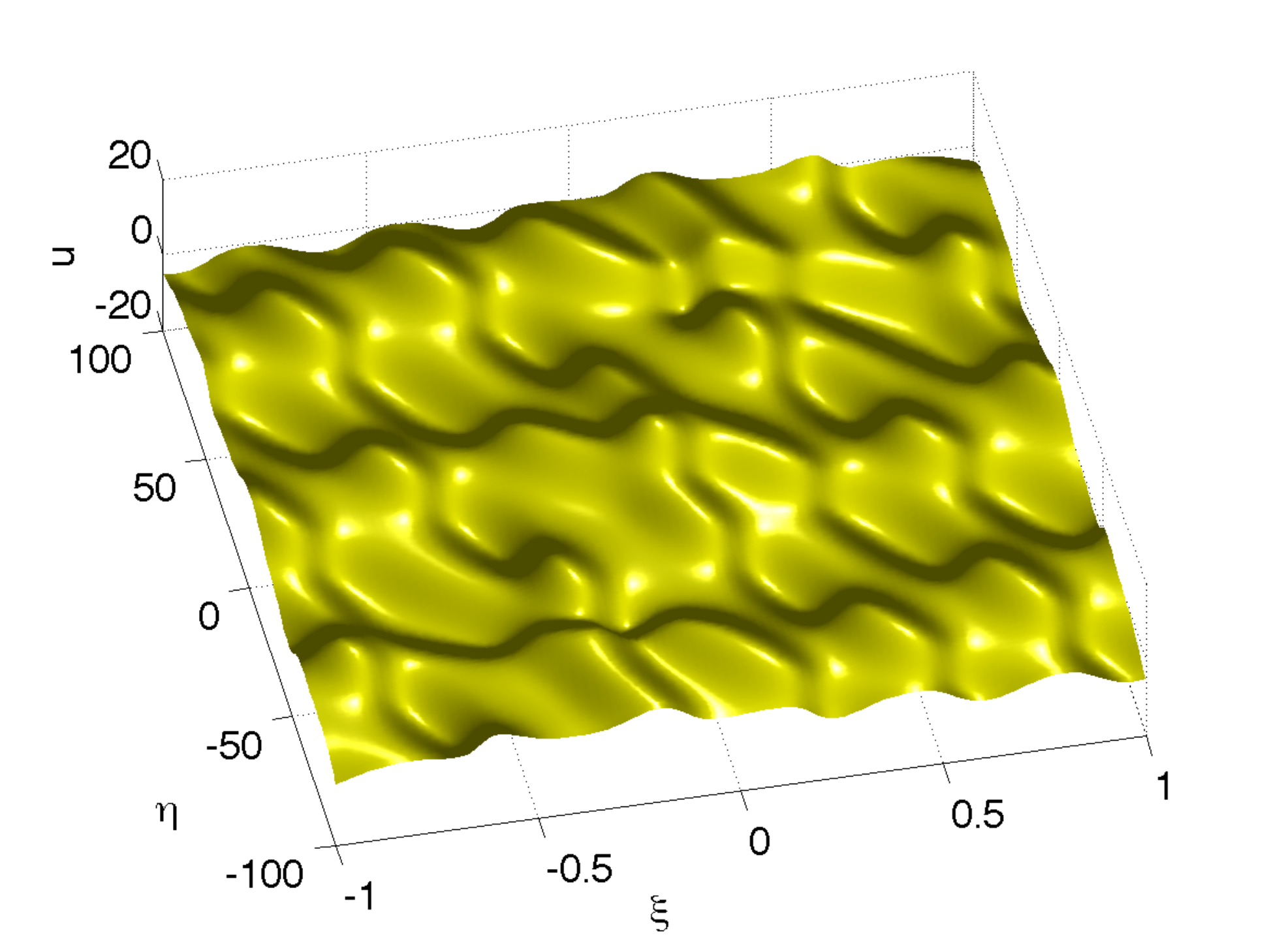}
  \includegraphics[width=0.49\linewidth]{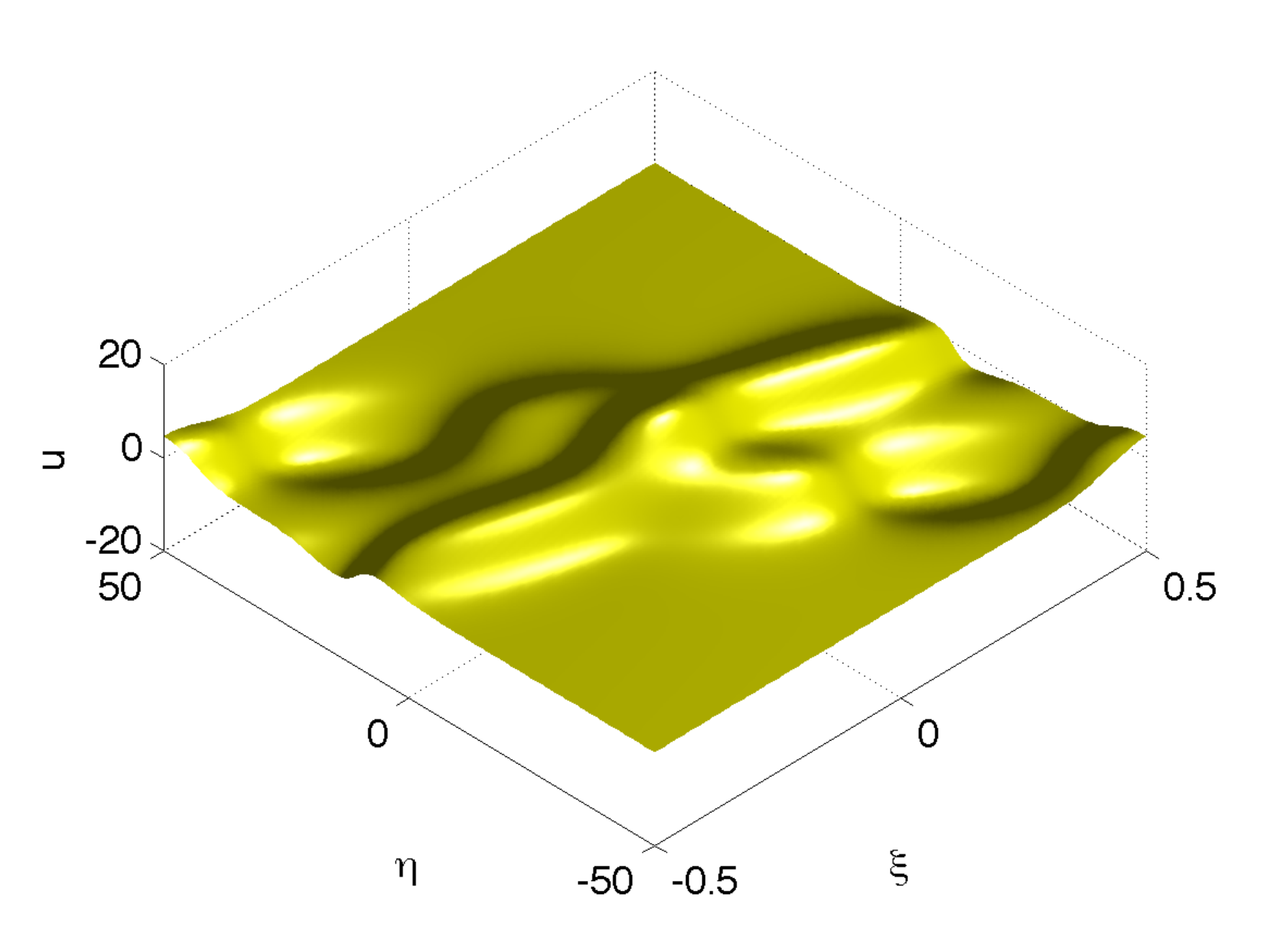}
  \caption{\label{SGimg4} Solution (\ref{SGformula}) to the SG equation (\ref{SG}) 
  on the surface with branch points (\ref{img4}); on the left for 
  $\epsilon=0.5$, on the right for 
  $\epsilon=10^{-12}$. }
\end{figure}

\section{Outlook}
In this paper, we have presented an algorithm to treat general 
hyperelliptic Riemann surfaces in Weierstrass form. The code is able 
to identify algorithmically a basis of the homology from a list of 
the branch points which was demonstrated for random points. It was shown that 
machine precision can be reached even in almost degenerate 
situations. This made it possible to study numerically the solitonic 
limit of solutions to the SG equation in terms of multi-dimensional 
theta functions.

The efficiency of the code allows the numerical study of functions on the 
modular space of  hyperelliptic surfaces, see \cite{Sarnak} for 
elliptic modular invariants. In \cite{klkoko}, extremal 
properties of the determinant of the Laplacian in the Bergman metric 
on the modular space of genus 2 Riemann surfaces were considered. This 
determinant could be given, however, in terms of theta functions 
only. The present code would allow the study of similar questions for 
modular invariants involving integrals over the whole Riemann surface 
as Faltings' $\delta$-invariant \cite{Faltings} or the invariant studied in 
\cite{Jong}. This will be the subject of further research.

\end{document}